\documentclass[twoside,a4paper]{article}

\usepackage{a4wide}
\usepackage{fancyhdr}
\usepackage[english]{babel}
\usepackage{amssymb}
\usepackage{amsmath}
\usepackage{eucal}
\usepackage{theorem}
\usepackage{epsfig}

\theorembodyfont{\itshape}
\theoremheaderfont{\scshape}
\newtheorem{theor}{Theorem}
\newtheorem{propo}{Proposition}[section]
\newtheorem{lemma}[propo]{Lemma}
\newenvironment{proof}{\noindent{\scshape Proof.}}{\hspace{2mm} $\square$ \\}
\newenvironment{sketch}{\noindent{\scshape Sketch of the proof.}}{\hspace{2mm} $\square$ \\}

\newcommand{\E}{\mathbf{E}}
\newcommand{\R}{\mathbf{R}}
\newcommand{\Z}{\mathbf{Z}}
\newcommand{\super}[1]{^{^{_{\,#1}}}}
\newcommand{\fpa}{first particle}
\newcommand{\dip}{distinguished particle}
\newcommand{\rpa}{red particle}
\newcommand{\gpa}{green particle}
\newcommand{\rep}{renewal point}
\newcommand{\stc}{strongly connected}
\newcommand{\wec}{weakly connected}

\newcommand{\phasea}   {1}
\newcommand{\phaseb}   {2}
\newcommand{\treefig}  {3}
\newcommand{\dualfig}  {4}
\newcommand{\renewfig} {5}
\newcommand{\betafig}  {6}
\newcommand{\fpafig}   {7}
\newcommand{\triangfig}{8}

\DeclareMathOperator{\Var}{Var}
\DeclareMathOperator{\card}{card}

\begin{document}
\thispagestyle{empty}

\begin{center}
\textbf{\large ECOLOGICAL SUCCESSION MODEL}
\footnote{\hspace{-16pt} \textit{American Mathematical Society} 1991 \textit{subject classifications}. 60K35} \\
\footnote{\hspace{-16pt} \textit{Key words and phrases}. Competition model, phase transition, interacting particle system.} \\ \vspace{8pt}
\textsc{By Nicolas Lanchier} \\ \vspace{4pt}
\textit{Universit\'e de Rouen} \\ \vspace{4pt}
\end{center}


\begin{abstract}
\noindent
 We introduce in this paper a new interacting particle system intended to model an example of ecological succession involving two species:
 \emph{Pteridium aquilinum L.} commonly called bracken and \emph{Fagus sylvatica L.} or european beech.
 The species both compete over vacant sites, assimilated to windfall, according to fundamentally different evolutionary strategies.
 Our main objective is to exhibit phase transitions for this process by proving that there exist three possible evolutions of the system to distinct
 ecological balances.
 More precisely, the whole population may die out, the european beech may conquer the bracken, and coexistence of both species may occur.
\end{abstract}


\section{\normalsize\bf{Introduction}}
\label{introduction}

\noindent
 The ecological model that we introduce in this article is a Markov process in which the state at time $t$ is given by some configuration
 $\eta_t : \Z^d \longrightarrow \{0, 1, 2 \}$.
 A site $x \in \Z^d$ is said to be empty or vacant if $\eta (x) = 0$, occupied by \emph{Pteridium aquilinum L.} or bracken if $\eta (x) = 1$,
 and occupied by \emph{Fagus sylvatica L.} or european beech if $\eta (x) = 2$.
 If the site $x$ is empty, it will be called sometimes a windfall.
 In what follows, we will formulate our results in the ecological context.
 However, to make the proofs more explicit, we will use another terminology, replacing the brackens by \emph{\rpa s} and the european beeches by
 \emph{\gpa s}.
 The evolutionary strategies of the two species involved in the competition are quite different and can be described as follows.

\begin{list}{}{\setlength{\leftmargin}{12pt} \setlength{\itemsep}{0pt} \setlength{\parsep}{2pt}}
 \item[1.] The brackens are rhizomatical ferns that colonize their neighboring sites thanks to an underground plexus of rhizomes.
  Moreover, their optimal development requires sunlight so that they can colonize the windfalls only.
 \item[2.] The european beeches produce an important amount of seeds that are spatially dispersed by animals (zoochory) and can remain
  in the soil as seed bank ; the development occurs on the shadow only and so, close to the adult trees or sheltered from the brackens.
\end{list}

\noindent
 With this in mind, we can now describe the dynamics of the ecological succession.
 Let $\mathcal V$ be the set of neighbors of 0 in $\Z^d$ given by $\mathcal V = \{\,y \in \Z^d \,,\,|\!|\,y \,|\!| \leq M \,\}$ for some positive
 integer $M$ and some norm $|\!| \cdot |\!|$ on $\R^d$.

\begin{list}{}{\setlength{\leftmargin}{12pt} \setlength{\itemsep}{0pt} \setlength{\parsep}{2pt}}
 \item[1.] A bracken (resp. a beech) standing at $x$ tries to give birth to its progeniture on the neighborood $x + \mathcal V$ at rate
  $\lambda_1$ (resp. $\lambda_2$).
  The new individual is sent to a site chosen at random from the neighboring sites.
 \item[2.] If the target site is a windfall, the birth occurs.
  Otherwise, it is suppressed.
 \item[3.] A bracken gives way to a new beech at rate 1.
 \item[4.] Finally, a beech lives for an exponentially distributed amount of time with mean 1 and then dies out to give way to a windfall.
\end{list}

\noindent
 For $i \in \{1, 2 \}$, denote by $v_i (x, \eta)$ the proportion of neighbors of $x$ occupied in the configuration $\eta$ by a type $i$ particle.
 Then the points 1 and 2 are equivalent to the fact that a windfall becomes occupied by a particle of type $i$ at rate $\lambda_i \,v_i (x, \eta)$. \\

\noindent
 We now formulate our results and construct step by step the phase space of the process in terms of the parameters $\lambda_1$ and $\lambda_2$.
 First of all, observe that if we start the evolution with a population of european beeches in the absence of brackens, the process
 we obtain is called the basic contact process with parameter $\lambda_2$.
 We know that, in this case, there exists a critical value $\lambda_c \in (0, \infty)$ such as the following holds.
 If $\lambda_2 \leq \lambda_c$ then the process converges in distribution to the all empty state, else there exists a stationary measure that
 concentrates on configurations with infinitely many beeches.
 See for instance Liggett (1999), Sect. 1.2. \\

\noindent
 Intuitively, we can think that, on the set where coexistence of both species does not occur, this result goes on holding for the
 ecological succession model.
 If $\lambda_2 \leq \lambda_c$, it is obvious, i.e. both species die out or coexist.
 If we now suppose that $\lambda_1 \geq \lambda_2 > \lambda_c$, we can prove that starting with infinitely many occupied sites,
 the process converges in distribution to a stationary measure that concentrates on configurations with infinitely many beeches.
 To do this, consider two initial configurations $\eta\super{1}$ and $\eta\super{2}$ such as $\eta\super{2} (x) = 0$ if $\eta\super{1} (x) = 0$,
 and $\eta\super{2} (x) = 2$ otherwise, i.e. $\eta\super{2}$ can be deduced from $\eta\super{1}$ by replacing the brackens by beeches.
 Observe that this implies that $\eta_t\super{2}$ is the basic contact process with parameter $\lambda_2$.
 Then by using a standard argument due to Harris (1972) we may run both processes simultaneously in such a way that if $\eta_t\super{2} (x) = 2$
 then $\eta_t\super{1} (x) \neq 0$.
 Hence, we can finally state the following

\begin{theor}
\label{dichotomie}
 Starting with infinitely many brackens, the process $\eta_t$ converges in distribution
\vspace{-8pt}
\begin{enumerate}
 \item[a.] either to the all empty state or to a measure that concentrates on configurations with infinitely many brackens (and beeches)
  if $\lambda_2 \leq \lambda_c$, \vspace{-8pt}
 \item[b.] either to the upper invariant measure of the supercritical contact process with parameter $\lambda_2$ or to a measure that concentrates on
  configurations with infinitely many brackens (and beeches) if $\lambda_1 \geq \lambda_2 > \lambda_c$.
\end{enumerate}
\end{theor}

\noindent
 Nevertheless, such a coupling fails to prove an analogous result in the case $\lambda_1 < \lambda_2$.
 But, since the particular dynamics of the ecological succession is inclined to favour the development of the beeches, it looks reasonable
 to think that, in the case $\lambda_2 > \lambda_c$, the beeches take over provided that $\lambda_1 \leq \lambda_2$.
 In fact, relying on some technicalities introduced by Neuhauser (1992), we can prove the two following theorems.

\begin{theor}
\label{diagonale}
 We suppose that $\lambda_2 > \lambda_c$ and $\eta_0$ is translation invariant.
\vspace{-8pt}
\begin{enumerate}
\item[a.] If $\lambda_1 < \lambda_2$ then, starting the evolution with infinitely many beeches, the brackens die out,
 i.e. $\eta_t \Rightarrow \mu_2$ the upper invariant measure of the contact process with parameter $\lambda_2$. \vspace{-8pt}
\item[b.] If $\lambda_1 = \lambda_2$ then, starting the evolution with infinitely many occupied sites, $\eta_t \Rightarrow \mu_2$.
 In particular, the beeches take over.
\end{enumerate}
\end{theor}

\noindent
 The last step is to show that the ecological succession does not reduce to the basic contact process, i.e. that coexistence of both species occurs for
 an open set of the parameters $\lambda_1$ and $\lambda_2$.
 If $\lambda_2 = 0$ and $d = 2$, the process is known as the forest fire model.
 See for instance Durrett and Neuhauser (1991).
 In this context, 0 = alive, 1 = on fire and 2 = burnt.
 We start the evolution with one burning tree in a virgin pine forest and follow the spread of the fire.
 A burning tree sends out sparks at rate $\lambda_1 = \alpha$ to one of its nearest neighbors $x + e_1$, $x + e_2$, $x - e_1$ or $x - e_2$ and burns for
 an exponentially distributed amount of time with parameter 1.
 Finally, a new pine regrows spontaneously on the burnt site at rate $\beta > 0$.
 By comparing this particle system with a one-dependent oriented percolation process, Durrett and Neuhauser (1991) have proved that, for fixed
 $\beta > 0$, there exists a critical value $\alpha_c \in (0, \infty)$ such as the fire has positive probability of not going out
 if $\alpha > \alpha_c$, i.e. there exists a nontrivial stationary distribution.
 Relying on their result and using an idea of Schinazi (2001) we can prove that this goes on holding for the ecological model with small values
 of $\lambda_2 > 0$.
 More precisely, we have the following

\begin{theor}
\label{haut}
 We suppose that $d = 2$ and that $\mathcal V$ reduces to the nearest neighbors.
 If $\lambda_1 > \alpha_c$ there exists a critical value $\beta_c \in (0, \infty)$ such as coexistence occurs as soon as $\lambda_2 \leq \beta_c$,
 i.e. there exists a stationary distribution that concentrates on configurations with infinitely many brackens.
\end{theor}

\pagebreak

\noindent
 In conclusion, there exist three possible evolutions of the system to distinct ecological balances.
 That is, depending on the parameters $\lambda_1$ and $\lambda_2$, the whole population may die out, the beeches may win the competition, and
 coexistence of both species may occur.
 For an picture of the space phase, see Fig. \phasea.
 On the other hand, it looks reasonable to conjecture, in view of the ecological dynamics, that for any $\lambda_1 > \alpha_c$
\begin{enumerate}
 \item There exists a critical value $\beta_c \,(\lambda_1) \in (0, \infty)$ such as conclusion of Theorem \ref{haut} holds if
       and only if $\lambda_2 \leq \beta_c \,(\lambda_1)$.
 \item The critical value $\beta_c \,(\lambda_1)$ is strictly increasing with respect to $\lambda_1$.
\end{enumerate}

\vspace{15pt}
\begin{figure}[ht]
\centering
\scalebox{0.4}{\input{phase-1.pstex_t}}
\end{figure}
\vspace{15pt}

\noindent
 Putting things together, the phase space of the particle system can be divided into three connex components corresponding to the three possible
 evolutions of the process.
 See Fig. \phaseb \ for a picture.
 This also implies that there exists a critical interval $I_c$ such as for any $\lambda_1 \in I_c$, the three possible evolutions may occur.
 In particular, we can notice that for any $\lambda_1 \in I_c$, an increase of the beeches birth rate $\lambda_2$ can make disappear the brackens and
 the beeches themselves.
 In the case $d = 1$ and $\mathcal V = \{-1, 1\}$, we have observed this phenomenon by simulation for any $\lambda_1 > \alpha_c$
 with $\alpha_c \simeq 7,95 \times \card \mathcal V = 15,9$.
 In particular, in view of the special geometry of the one dimensional process, we think that coexistence cannot occur in $d = 1$
 whenever $\lambda_2$ is supercritical. \\

\noindent
 The article is organized as follows:
 In Sect. \ref{construction}, we describe in greater details the graphical construction of the dual process using an idea of Harris (1972).
 In Sect. \ref{renewal}, we prove Theorem \ref{diagonale} for $\lambda_1 < \lambda_2$ in any dimension relying on the graphical construction given
 above and some results of Neuhauser (1992).
 The proof of Theorem \ref{diagonale} for $\lambda_1 = \lambda_2$, essentially based on classical properties of random walks on $\Z^d$, is carried
 out in two steps.
 In view of the well-known specificities of symetric random walks depending on the dimension, we start, in Sect. \ref{recurrence}, by proving
 the theorem in $d \leq 2$, and conclude, in Sect. \ref{transience}, dealing with the case $d \geq 3$.
 Finally, Sect. \ref{coexistence} is devoted to the description of results of Durrett and Neuhauser (1991) and the proof of Theorem \ref{haut}.

\vspace{15pt}
\begin{figure}[ht]
\centering
\scalebox{0.4}{\input{phase-2.pstex_t}}
\end{figure}
\vspace{15pt}


\section{\normalsize\bf{Construction and properties of the dual process}}
\label{construction}

\noindent
 We begin by constructing the process from a collection of Poisson processes in the case $\lambda_1 \leq \lambda_2$.
 For $x$, $y \in \Z^d$, $x - y \in \mathcal V$, let $\{\, T_n^{x,y} \,;\, n \geq 1 \,\}$ and $\{\, U_n^x \,;\, n \geq 1 \,\}$ be the arrival times
 of Poisson processes with rates $\lambda_2 \cdot (\card \mathcal V)^{-1}$ and 1 respectively.
 At times $T_n^{x,y}$, we draw an arrow from $x$ to $y$, toss a coin with success probability $(\lambda_2 - \lambda_1) / \lambda_2$, and, if there is a
 success, label the arrow with a 2.
 If $x$ is occupied, $y$ is vacant and the arrow is unlabelled, then a birth will occur in $y$.
 If the arrow is labelled with a 2, then the beeches only can give birth through this arrow.
 Finally, at times $U_n^x$, we put a $\times$ at $x$.
 In view of the dynamics, the effect of a $\times$ is to kill beeches and to permute the brackens into beeches.
 A result of Harris (1972) implies that such a graphical representation can be used to construct the process starting from any
 $\eta_0 \in \{0, 1, 2 \}^{\Z^d}$. \\

\noindent
 We say that $(x,0)$ and $(y,t)$ are \emph{\stc}, which we note $(x,0) \rightarrow (y,t)$, if there is a sequence of times
 $s_0 = 0 < s_1 < s_2 < \cdots < s_n < s_{n + 1} = t$ and spatial locations $x_0 \, = \,x, \ x_1, \ x_2, \ \cdots \ x_n \, = \,y$ so that
\begin{enumerate}
\item [(1)] for $1 \leq i \leq n$, there is an arrow from $x_{i - 1}$ to $x_i$ at time $s_i$ and
\item [(2)] for $0 \leq i \leq n$, the vertical segments $\{x_i\} \times (s_i, s_{i + 1})$ do not contain any $\times$.
\end{enumerate}
 If instead of (2) we have the condition
\begin{enumerate}
\item [(3)] the set $\displaystyle \bigcup_{i = 0}^n \ \{x_i\} \times (s_i, s_{i + 1})$ contains exactly one $\times$,
\end{enumerate}
 we say that $(x,0)$ and $(y,t)$ are \emph{\wec}, which we note $(x,0) \rightharpoonup (y,t)$. \\

\noindent
 After constructing the graphical representation we now define the dual process.
 Since the $\times$ do not kill the brackens but convert them into beeches, we must take into account the paths that contain at most one $\times$.
 So, to construct the dual process, we reverse the arrows and time by mapping $\tilde s = t - s$, and let
 $$ \tilde \eta_t^x \ = \ \{\, y \in \Z^d \,;\, (x,\tilde 0) \rightarrow (y,\tilde t) \ \textrm{or} \ (x,\tilde 0) \rightharpoonup (y,\tilde t) \,\}. $$
 Since it is in general easier to work with a forward process, we replace $\tilde \eta_t^x$ by the dual process $\hat \eta_t^x$ that is constructed by
 restoring the direction of time and arrows and by letting 
 $$ \hat \eta_t^x \ = \ \{\, y \in \Z^d \,;\, (x,0) \rightarrow (y,t) \ \textrm{or} \ (x,0) \rightharpoonup (y,t) \,\}. $$
 As in the basic contact process, $\tilde \eta_t^x$ and $\hat \eta_t^x$ have the same distribution. \\

\noindent
 Before sketching the proof of Theorem \ref{diagonale}, which is based on the graphical representation, we give some definitions intended to describe the
 geometry of the dual process.
 First of all, note that $\{\,\hat \eta_s^x \,;\, 0 \leq s < t \,\}$ has a tree structure composed of the points that are either strongly or weakly
 connected with $(x, 0)$.
 To figure on this distinction, denote by $\Gamma$ the tree formed by the set of points in $\Z^d \times \R^+$ that are \stc \ with $(x, 0)$, i.e.
 $$ \Gamma \ = \ \{\, (y, t) \in \Z^d \times \R^+ \,;\, (x, 0) \rightarrow (y, t) \,\}. $$
 Clearly, the growth of $\Gamma$ is broken at some points by a $\times$ from which a new tree takes form so that we obtain the picture of a pyramid of
 trees connected with $\Gamma$ by some $\times$.
 In the following, $\Gamma$ will be called the \emph{upper tree starting at} $(x, 0)$ and the trees starting at a $\times$ the \emph{lower trees}.
 For an illustration, see Fig. \treefig.
 To motivate such a breaking up, note that, in view of the translation invariance of the graphical representation in space and time, the upper tree and
 the lower trees are identically distributed.
 Furthermore, $\Gamma$ has obviously the same distribution as the tree structure of the contact process with parameter $\lambda_2$, that is the set
 of points in $\Z^d \times \R^+$ that can be connected by some path to $(x, 0)$.
 See for instance Durrett (1995), Sect. 3 or Liggett (1999), Sect. 1.1 for a construction of the basic contact process.
 So, one of the main interest of our decomposition is that it allows to break up the tree starting at $(x, 0)$ into identically distributed pieces
 that have some usefull similarities with some well-known structures. \\

\noindent
 Now, denote by $\hat \eta_t^{x,1}$ and $\hat \eta_t^{x,2}$ the subsets of the dual process given by
 $$ \hat \eta_t^{x,1} \ = \ \{\, y \in \Z^d \,;\, (x,0) \rightarrow (y,t) \,\} \qquad \textrm{and} \qquad
    \hat \eta_t^{x,2} \ = \ \{\, y \in \Z^d \,;\, (x,0) \rightharpoonup (y,t) \,\}. $$
 By analogy with the breaking up of the tree structure in upper tree and lower trees, the elements of $\hat \eta_t^{x,1}$ and $\hat \eta_t^{x,2}$ will be
 respectively called \emph{upper ancestors} and \emph{lower ancestors}.
 As in the multitype contact process, the tree structure of the dual process $\{\,\hat \eta_s^x \,;\, 0 \leq s < t \,\}$ allows to define an ancestor
 hierarchy in which the members are arranged according to the order they determine the color of $(x, 0)$.
 Here, the special geometry of the dual described above plays a leading role since the color of $(x, 0)$ also depends on the nature of the ancestors.
 To specify this idea, we now describe in greater detail the hierarchy. \\

\noindent
 First of all, suppose that $\lambda_1 = \lambda_2$.
 Denote by $\hat \eta_t^x (n)$ the $n$-th member of the ordered ancestor set and let $\hat \eta_t^{x,1} (k) = \hat \eta_t^x (n_k)$ be the $k$-th
 upper ancestor.
 Later on, $\hat \eta_t^{x,1} (1)$ will be called the \emph{\dip}, and the first ancestor the \emph{\fpa}.
 If the \fpa \ is a upper ancestor that lands on a \rpa \ (resp. a \gpa), then it will paint $(x, 0)$ in red (resp. green).
 If it is a lower ancestor that lands on a red site, in view of the passage through one $\times$, it will bring a \gpa \ to $(x, 0)$.
 If it is a lower ancestor that lands on a green site, the path taken by the \gpa \ until it is destroyed by a $\times$ can prevent the rising of some
 ancestors toward $(x, 0)$.
 Hence, we discard all the ancestors whose path crosses the \gpa.
 Then, we repeat the same procedure with next remaining ancestor, and so on. \\

\begin{figure}[h]
\centering
\scalebox{0.4}{\input{uptree.pstex_t}} \hspace{50pt}
\scalebox{0.4}{\input{dual.pstex_t}}
\end{figure}

\noindent
 If $\lambda_2 > \lambda_1$, the arrows labelled with a 2 are forbidden for the \rpa s.
 So, if $\hat \eta_t^x (1)$ does not cross any 2-arrows or lands on a green site or a windfall, we apply the same procedure as before.
 If it is a lower ancestor that lands on a red site and that crosses a $\times$ before the first 2-arrow, it will paint $(x, 0)$ in green.
 In the other cases, we follow the path taken by the particle until we first cross a 2-arrow.
 Then we discard all the ancestors of the point where this arrow is attached, and start afresh with the first remaining ancestor, and so on.
 For instance, in Fig. \dualfig, the \fpa \ and the third ancestor are lower ancestors that land on a green site so they cannot paint $(x, 0)$ any color.
 Moreover, the path of the second ancestor is forbidden by the rising of a \gpa \ at site $x + 1$ so that the \dip \ $\hat \eta_t^x (4)$
 paints $(x, 0)$ in red.
 If we now suppose that $\lambda_2 > \lambda_1$, we can notice that $\hat \eta_t^x (4)$ also fails and that the fifth ancestor paints $(x, 0)$ in green.


\section{\normalsize\bf{Proof of Theorem \ref{diagonale} for $\lambda_1 < \lambda_2$}}
\label{renewal}

\noindent
 In this section, we will prove Theorem \ref{diagonale} for $\lambda_1 < \lambda_2$ in any dimension.
 The strategy of the proof is quite simple ; it relies on a usefull idea of Neuhauser (1992) that consists in breaking up the path of the \dip \
 at certain points into i.i.d. pieces.
 We will first remind this result, describing the intuitive idea on which it is based, and then conclude that the brackens die out as soon
 as $\lambda_1 < \lambda_2$. \\

\noindent
 To begin with, observe that if the upper tree starting at $(x, 0)$ does not live forever then, for $t$ large enough, $\hat \eta_t^x = \hat \eta_t^{x,2}$.
 In particular, in view of the passage of any lower ancestor through one $\times$, $(x, 0)$ cannot be reached by a \rpa \ (see the description of the
 ancestor hierarchy in Sect. \ref{construction}).
 So we can suppose from now on that the upper tree $\Gamma$ lives forever.
 Note that the probability of such an event is given by the survival probability of the contact process with parameter $\lambda_2$ starting from one
 infected site, that is obviously positive since $\lambda_2$ is supercritical.
 In this case, we proceed in two steps by firstly checking that the path of the \dip \ is forbidden for the red, and then proving that we can bring
 a \gpa \ to $(x, 0)$.
 To do this, we define an embedded random walk for the \dip \ by breaking its evolution at some points called \emph{\rep s}. \\

\noindent
 We now specify this construction.
 First of all, denote by $s_k$, $k \geq 1$, the jumping times of the \dip.
 Then, follow the path of the particle inside the upper tree and, whenever it jumps to a new site
 $\hat \eta_{s_k}^{x,1} (1)$, look at the branch starting at $(\hat \eta_{s_k}^{x,1} (1), s_k)$.
 If this branch lives forever, $(\hat \eta_{s_k}^{x,1} (1), s_k)$ will be called a \rep.
 Let $(S_n \super{0}, T_n \super{0})$ be the location of the $n$-th renewal.
 For an illustration, see Fig. \renewfig.
 Denote by $X_i \super{0}$ the spatial displacement between consecutive \rep s, and by $\tau_i \super{0}$ the corresponding temporal displacement so that
 $$ S_n \super{0} \ = \ x \ + \ \sum_{i = 1}^n \ X_i \super{0}  \qquad \textrm{and} \qquad  T_n \super{0} \ = \ \sum_{i = 1}^n \ \tau_i \super{0}. $$
 Our main ingredient to prove Theorem \ref{diagonale} is given by the following proposition.
 We just describe the intuition behind this result.
 For the details of the proof, see Neuhauser (1992), Sect. 2.

\begin{propo}
\label{random-walk}
 If the upper tree lives forever, $\{(X_i \super{0}, \tau_i \super{0})\}_{i \geq 1}$ form an i.i.d. family of random \linebreak vectors on $\Z^d \times \R^+$.
 Moreover, the tail distributions of $X_i \super{0}$ and $\tau_i \super{0}$ have exponential bounds, i.e. there exist positive constants $C$ and $\beta$
 such as
 $$ P \,(\,|X_i \super{0}| > t \,) \ \leq \ C \,e^{- \beta t}  \qquad \textrm{and} \qquad  P \,(\,\tau_i \super{0} > t \,) \ \leq \ C \,e^{- \beta t} $$
 holds for all $t \geq 0$.
\end{propo}

\begin{sketch}
 First of all, denote by $\sigma_0$ the first jumping time of the \dip, that is
 $\sigma_0 = \inf \,\{\,t > 0 \,; \,\hat \eta_t^{x,1} (1) \textrm{ hits a } \times \,\}$, by $x_1$ its spatial location after $\sigma_0$, and by $\beta_1$
 the branch of $\Gamma$ starting at $(x_1, \sigma_0)$.
 For a picture, see Fig. \betafig.
 If $\beta_1$ lives forever then $(x_1, \sigma_0)$ is the first renewal point.
 Else, we define the sequences $\{(x_k, \sigma_{k - 1})\}_{k \geq 1}$ and $(\beta_k)_{k \geq 1}$ as follows.
 If $\sigma_{k - 1} < \infty$, let $x_k$ be the location of the particle after $\sigma_{k - 1}$ and $\beta_k$ the branch starting
 at $(x_k, \sigma_{k - 1})$.
 Note that such a branch always exists since we have supposed that $\Gamma$ lives forever.
 Then denote by $\sigma_k$ the time when $\beta_k$ dies out.
 The sequences are defined until $\sigma_k$ is equal to infinity.
 If $\sigma_k = \infty$ then $\beta_k$ lives forever and $(x_k, \sigma_{k - 1})$ is the first renewal point of the \dip.
 To find the next one, we start over again the whole procedure replacing $(x, 0)$ by $(x_k, \sigma_{k - 1})$, and so on.
 Now, look at the path of the particle between time 0 and time $\sigma_{k - 1}$.
 Whenever the particle jumps within a new branch $\beta_i$, the location of its target site $x_i$ clearly depends on the past history of $\Gamma$,
 that is on the geometry of $\Gamma$ between 0 and $\sigma_{i - 1}$.
 On the contrary, from time $\sigma_{k - 1}$, the particle stays forever inside $\beta_k$ and its path only depends on this branch so that what
 happens before and after a renewal point is determined by disjoint parts of $\Gamma$.
 This implies that the random vectors $(X_i \super{0}, \tau_i \super{0})$ are independent.
 Finally, since the graphical gadget is translation invariant in space and time the vectors $(X_i \super{0}, \tau_i \super{0})$ are also identically
 distributed.
 The proof of the exponential bounds, more technical, can be found in Neuhauser (1992), Sect. 2.
\end{sketch}

\noindent
 By Proposition \ref{random-walk}, the sequence $(S_n \super{0}, T_n \super{0})$ defines an embedded random walk for the particle.
 Moreover, the exponential bounds on the displacements $X_i \super{0}$ and $\tau_i \super{0}$ give us control over the location of the particle between
 consecutive \rep s.
 Since the branch starting at each renewal grows at most linearly in time, we know that the particle stays within a \emph{triangle} whose base
 contains the following renewal. For a picture, see Fig. \triangfig. \\

\pagebreak

\noindent
 A triangle is said \emph{closed for the red} if the \dip \ goes through an arrow labelled with a 2 to cross it.
 Since $\lambda_1 < \lambda_2$, the event that a triangle is closed for the red obviously occurs with positive probability.
 So, by the Borel-Cantelli Lemma, choosing $t$ large enough, we can find arbitrarily many closed triangles. \\

\begin{figure}[ht]
\centering
\scalebox{0.4}{\input{renew.pstex_t}} \hspace{50pt}
\scalebox{0.4}{\input{beta.pstex_t}}
\end{figure}

\noindent
 To conclude the proof, we now use the dual process $\tilde \eta_s\super{(x, t)}$, $0 \leq s \leq t$, starting at $(x, t)$ and determine the ancestor
 hierarchy after $t$ units of time by going backwards in time.
 We construct by induction a sequence of upper ancestors $\tilde \xi_t^{x,1} (k)$, $k \geq 1$, that are candidates to bring a \gpa \ to $(x, t)$.
 The first member of the sequence is the \dip \ $\tilde \eta_t^{x,1} (1)$.
 Then, we follow the path the \dip \ takes to paint $(x, t)$ by going forward in time until we first cross a 2-arrow.
 Finally, we look backwards in time starting from the location where this arrow is attached and discard all the offspring of this point (that is
 the next few members of the ordered ancestor set) and all the ancestors that land at time 0 on the same site as the \dip.
 The second member of the sequence is then the first upper ancestor that is left.
 We repeat the procedure with this ancestor, and so on until we run out of upper ancestors.
 Now, denote by $\xi_t$ the set of members of the sequence and by $B_s\super{(x, t)}$ the set of sites occupied at time $s$ by a \gpa.
 Since the tree $\Gamma$ is linearly growing in time, we can make the cardinality of $\xi_t$ arbitrarily large by choosing $t$ large enough.
 In particular, given $\varepsilon > 0$ and $M > 0$ there exists a time $t_0 \geq 0$ so that
 $P \,(\,\card \,(\xi_t) < M \,) \leq \varepsilon$ for any $t \geq t_0$.
 By Lemma 9.14 of Harris (1976), this together with the translation invariance of $\eta_0$ and the fact that $\eta_0 (x) = 2$ with positive
 probability, imply that
 $$ \lim_{t \rightarrow \infty} \ P \,(\,\xi_{t - 1} \,\cap \,B_1\super{(x, t)} = \,\varnothing \,) \ = \ 0. $$
 Hence, for $t$ large enough, at least one candidate lands on a green site.
 Now, consider the first member of the sequence $\tilde \xi_t^{x,1} (k)$ that lands on a site occupied by a \gpa \ and follow the path of this
 particle by going forward in time.
 Then, whenever it crosses an arrow, its target site is either empty or already occupied by a \gpa \ (resulting from the passage of a \rpa \ through
 the first $\times$ located under the arrowhead).
 Thus, $(x, t)$ will be eventually painted in green.
 This completes the proof of Theorem \ref{diagonale} for $\lambda_1 < \lambda_2$.


\section{\normalsize\bf{Proof of Theorem \ref{diagonale} for $\lambda_1 = \lambda_2$ in dimension $\leq 2$}}
\label{recurrence}

\noindent
 The technicalities and tools we will make use to prove Theorem \ref{diagonale} in the case $\lambda_1 = \lambda_2$ are quite different depending on the
 dimension of the state space.
 In this section, we will deal with the case $d \leq 2$, relying on the recurrence of 1 and 2 dimensional random walks.
 To begin with, we will show that the \fpa \ can be trapped with probability one inside a lower tree that lives forever, so that, for $t$ large enough,
 the first ancestor is a lower ancestor.
 At this point, the worst scenario we have in mind is that the \dip \ lands on a red site and the first ancestor on a green one.
 In such a case the \fpa, which is a lower ancestor, cannot paint $(x, 0)$ any color whereas the \dip \ can possibly bring a \rpa \ to $(x, 0)$.
 To conclude, we will then prove that this bad event is negligible showing that, with probability one, we can make coalesce the \dip \ and the \fpa \
 and so make them land on the same site.
 If they both land on a green (resp. on a red), the \dip \ (resp. the \fpa) will paint $(x, 0)$ in green, and the theorem will follow. \\

\noindent
 From now on, we suppose that $\lambda_1 = \lambda_2$ and denote by $\lambda$ their common value.
 The first step is given by the lemma

\begin{lemma}
\label{fpa}
 There exists a time $\theta_1$ a.s. finite such as for any $t \geq \theta_1$, the \fpa \ is a lower ancestor.
\end{lemma}

\begin{proof}
 To begin with, let $s_k$, $k \geq 1$, be the jumping times of the \fpa \ and $x_k$ its location before $s_k$.
 Denote by $\sigma_1$ the first time the particle crosses a $\times$, i.e.
 $$ \sigma_1 \ = \ \inf \,\{\, t \geq 0 \,; \,\hat \eta_t^x (1) \ \textrm{is a lower ancestor} \,\}, $$
 and by $\Omega_1$ the lower tree starting at $(x_1, \sigma_1)$, that is the first lower tree the particle visits.
 See Fig. \fpafig \ for an illustration.
 In view of the tree structure, once the particle penetrates in $\Omega_1$, it remains trapped inside (as lons as the lower tree is alive).
 Hence, if $\Omega_1$ lives forever, $\hat \eta_t^x (1)$ is a lower ancestor for any $t \geq \sigma_1$, and the proof is done.
 Else, denote by $\sigma_2$ the first time the \fpa \ visits a new lower tree after $\Omega_1$ dies and by $\Omega_2$ the lower tree starting at
 $(\hat \eta_{\sigma_2}^x (1), \sigma_2)$.
 Note that for all $k \geq 1$ the path the particle takes to climb from $(x_k, s_k^-)$ to $(x, 0)$ contains one $\times$ so $\sigma_2$ is finite a.s.
 and $\Omega_2$ is well defined.
 While the particle is not trapped in a lower tree that lives forever, we thus construct by induction a sequence of trees $\Omega_k$ visited by
 $\hat \eta_t^x (1)$.
 Now, denote by $B_n$ the event that the first $n$ trees $\Omega_1$, $\Omega_2$, \ldots, $\Omega_n$ are bounded and, for any $k \geq 1$, by $A_k$ the
 event that the $k$-th tree lives forever.
 If $A_k$ does not occur then $\Omega_{k + 1}$ is well defined and the event $A_{k + 1}$ is determined by parts of the graph that are after $\Omega_k$
 dies so $A_k$ and $A_{k + 1}$ are independent.
 More generally, since the trees $\Omega_1$, $\Omega_2$, \ldots, $\Omega_{k + 1}$ are disjoint, $A_1$, $A_2$, \ldots, $A_{k + 1}$ are independent.
 Moreover, $\Omega_k$ has the same distribution as the tree structure of the contact process so the probability that $A_k$ occurs is given by
 $\rho_{\lambda}$, the survival probability of the contact process with parameter $\lambda$ starting from one infected site.
 This implies that
 $$ P \,(B_n) \ = \ P \,(\,A_1^c \,\cap \,\ldots \,\cap \,A_{n - 1}^c \,\cap \,A_n^c \,) \ = \
    \prod_{k = 1}^n \,P \,(\, A_k^c \,) \ = \ (1 - \rho_{\lambda})^n. $$
 Finally, since $\lambda$ is supercritical, the survival probability of the contact process $\rho_{\lambda}$ is strictly positive so that
 $\displaystyle \lim_{n \rightarrow \infty} P \,(B_n) = 0$.
 This completes the proof of the lemma.
\end{proof}

\noindent
 The next step of the proof is to show that the particles coalesce with probability one.
 To do this, note that after penetrating in a lower tree $\Gamma_1$ that never dies, the \fpa \ is \wec \ with $(x, 0)$ and so jumps within a new branch each
 time it meets a $\times$.
 In particular, from time $\theta_1$, the path of the particle can be broken up, as for the \dip, into i.i.d. pieces.
 We define the \rep s of the \fpa \ as before replacing the upper tree $\Gamma$ by the lower tree $\Gamma_1$ starting at $(S_0 \super{1}, \theta_1)$.
 Here, $S_0 \super{1}$ is the site where the particle jumps at time $\theta_1$.
 We denote by $(S_n \super{1}, T_n \super{1})$ the location of the $n$-th renewal after $\theta_1$, and by $X_i \super{1}$ and $\tau_i \super{1}$ the spatial
 and temporal displacements between two consecutive renewals, so
 $$ S_n \super{1} \ = \ S_0 \super{1} \ + \ \sum_{i = 1}^n \ X_i \super{1}  \qquad \textrm{and} \qquad
    T_n \super{1} \ = \ \theta_1 \ + \ \sum_{i = 1}^n \ \tau_i \super{1}. $$
 By translation invariance of the graphical representation (see the description of the tree structure in Sect. \ref{construction}), the families
 $\{(X_i \super{0}, \tau_i \super{0})\}_{i \geq 1}$ and $\{(X_i \super{1}, \tau_i \super{1})\}_{i \geq 1}$ are identically distributed so that
 Proposition \ref{random-walk} remains valid for the random vectors $(X_i \super{1}, \tau_i \super{1})$.
 In particular, as long as their triangles do not collide, both particles behave nearly like independent random walks.
 This constitutes the main ingredient to prove coalescence.
 We now specify this idea in greater details, describing the approach of Neuhauser (1992), Sect. 3. \\

\begin{figure}[ht]
\centering
\scalebox{0.4}{\input{fpa.pstex_t}} \hspace{60pt}
\scalebox{0.4}{\input{triang.pstex_t}}
\end{figure}

\noindent
 The first idea is to extend the notion of renewals for both particles, that is to break up the set of both paths into i.i.d. pieces.
 To do this, we say that an ancestor is \emph{good} at time $t$ if it did not meet any arrow since its last renewal.
 Observe that if both particles are good at the same time, what happens before and after that time uses disjoint parts of the graph and then is independent.
 We now prove that both particles are good i.o. at the same time.

\begin{lemma}
\label{good}
 $P \,(\,\textit{both particles are good at the same time i.o.}\,) = \, 1$.
\end{lemma}

\begin{proof}
 To begin with, we construct by induction two sequences of subscripts $(n_k)_{k \geq 0}$ and $(m_k)_{k \geq 1}$ as follows.
 We let $n_0 = 1$, and for any $k \geq 1$
 $$ m_k \ = \ \min \,\{\,m \geq 1 \,,\, T_m \super{0} > T_{n_{k - 1}} \super{1} \,\} \qquad \textrm{and} \qquad
    n_k \ = \ \min \,\{\,n \geq 1 \,,\, T_n \super{1} > T_{m_k} \super{0} \,\}. $$
 Then, denote by $A_k$ the event that the \fpa \ lives without giving birth between time $T_{n_{k - 1}} \super{1}$ and $T_{m_k} \super{0}$.
 Note that if $A_k$ occurs then both particles are obviously good at time $T_{m_k} \super{0}$.
 Moreover, in view of the exponential bound given by Proposition \ref{random-walk}, for any $T > 0$
 $$ P \,(\,T_{m_k} \super{0} - T_{n_{k - 1}} \super{1} > T \,) \ \leq \ P \,(\,T_{m_k} \super{0} - T_{m_k - 1} \super{0} > T \,) \ = \
    P \,(\,\tau_1 \super{0} > T \,) \ \leq \ C \,e^{- \beta T} $$
 so that $P \,(A_k) \ \geq \ (\,1 - C \,e^{- \beta T} \,) \ e^{- (1 + \lambda) T}$.
 Here, $e^{- (1 + \lambda) T}$ is the probability that the \fpa \ lives without giving birth for $T$ units of time.
 Since this holds for all $T > 0$, there exists a constant $\varepsilon > 0$ such as $P (A_k) \geq \varepsilon$.
 On the other hand, the events $A_k$ are determined by disjoint parts of the graphical gadget so they are independent.
 Finally, by the Borel-Cantelli Lemma, we can conclude that
 $$ P \,(\,\textrm{both particles are good at the same time i.o.}\,) \ \geq \ P \,( \ \limsup_{k \rightarrow \infty} A_k \,) \ = \ 1. $$
 This proves the lemma.
\end{proof}

\noindent
 To make coalesce both particles together, we now proceed in two steps.
 First of all, relying on the recurrence of 1 and 2 dimensional random walks, Neuhauser has proved that with positive probability we can bring both
 particles within a finite distance $K$ without collision of their triangles.
 See Neuhauser (1992), Sect. 3, for a proof.
 Then, as soon as the particles are close enough to each other, we try to make them coalesce.
 More precisely, we have the following

\begin{lemma}
 Let $K$ be a positive integer and suppose that both particles are within a distance $K$ at some time $t \geq \theta_1$.
 Then the event $B$ that the particles coalesce has positive probability.
\end{lemma}

\begin{proof}
 To find a lower bound for $P (B)$, we first require the \dip \ to keep still for $3 dK$ units of time.
 Then we can make both particles coalesce in less than $dK$ steps by increasing or decreasing each coordinate of the \fpa.
 To estimate this event, note that the probability of having neither birth nor death between time $t$ and $t + 1$, a good oriented arrow between
 time $t + 1$ and $t + 2$, and a death between time $t + 2$ and $t + 3$ is given by
 $$ e^{-2} \ (1 - e^{-1}) \ e^{- \lambda} \ |\mathcal V|^{-1} \ (1 - e^{- \lambda}). $$
 Since it takes the \fpa \ at most $3dK$ units of time to reach the \dip, we can conclude that
 $$ P \,(B) \ \geq \ e^{- 3dK (1 + \lambda)} \ [\,e^{-2} \ (1 - e^{-1}) \ e^{- \lambda} \ |\mathcal V|^{-1} \ (1 - e^{- \lambda}) \,]^{dK} \ > \ 0 $$
 where $e^{- 3dK (1 + \lambda)}$ is the probability that the \dip \ survives without giving birth for $3dK$ units of time.
 This completes the proof.
\end{proof}

\noindent
 If we do not succeed in gluing the particles together, we use the restart argument given by Lemma \ref{good}, i.e. we wait until both particles are good at the
 same time and then start over again the whole procedure.
 Since the set of paths is broken into i.i.d. pieces, we can apply the Borel-Cantelli Lemma to conclude that coalescence eventually occurs. \\

\noindent
 We prove Theorem \ref{diagonale} in $d \leq 2$ by using, as before, the dual process $\tilde \eta_s\super{(x, t)}$, $0 \leq s \leq t$, starting at $(x, t)$.
 Since both particles coalesce a.s., we can suppose, taking $t$ large enough, that they land at time 0 on the same site.
 If it is a red site, the \fpa, that is a lower ancestor by Lemma \ref{fpa}, will paint $(x, t)$ in green.
 If both particles land on a green site, we have already proved in Sect. \ref{renewal} that the \dip \ can bring a \gpa \ up to $(x, t)$.
 Finally, if the target site is a windfall, we start over again with the second ancestor, and so on.
 Since the tree starting at the coalescence point is linearly growing in time, we eventually find, by Harris (1976), an ancestor that lands on
 an occupied site and paints $(x, t)$ in green.


\section{\normalsize\bf{Proof of Theorem \ref{diagonale} for $\lambda_1 = \lambda_2$ in dimension $\geq 3$}}
\label{transience}

\noindent
 The strategy of the proof to deal with the case $d \geq 3$ is quite different.
 To begin with, we will construct by induction an ordered set of ancestors, $\gamma_t^x (k)$, $k \geq 1$, that are candidates to paint $(x, 0)$ in green.
 Using Lemma \ref{fpa}, we will prove that for any $k \geq 1$, and for $t$ large enough, $\gamma_t^x (k)$ is a lower ancestor that comes before the \dip \ in
 the ancestor hierarchy.
 Then, relying on the transience of $d$-dimensional random walks for $d \geq 3$, we will extract a subsequence of particles $\gamma_t^x (k_i)$, $i \geq 1$,
 that never coalesce.
 In particular, the number of sites occupied by the particles can be made arbitrarily large so that we can eventually find a particle
 landing on a red site that will paint $(x, 0)$ in green. \\

\noindent
 We start by constructing inductively the particle system $\gamma_t^x (k)$, $k \geq 1$.
 The first member of the sequence is the \fpa.
 Then, follow its path until it penetrates in a lower tree $\Gamma_1$ that lives forever.
 By Lemma \ref{fpa}, we know that this occurs at some time $\theta_1$ a.s. finite.
 Now, consider the ancestor hierarchy at that time and discard all the ancestors that land on $S_0 \super{1}$ or on a branch that does not live forever.
 We remind that $S_0 \super{1}$ is the spatial location of the \fpa \ at time $\theta_1$.
 The second member of the sequence is then the first ancestor that is left.
 Observe that such an ancestor obviously exists since we have supposed that $\Gamma$ lives forever.
 Moreover, using the arguments of \ref{fpa}, we can prove that it penetrates a.s. in a lower tree $\Gamma_2$ that never dies.
 We then repeat the same procedure to define the third particle, and so on.
 For any $k \geq 1$, let $\Gamma_k$ be the infinite lower tree visited by the $k$-th particle, $\theta_k$ the first time the particle penetrates in $\Gamma_k$
 and $S_0 \super{k}$ the spatial location of the particle at time $\theta_k$.
 As for the \fpa, we can break up the path of the $k$-th particle inside $\Gamma_k$ into i.i.d. pieces.
 We denote by $(S_n \super{k}, T_n \super{k})$ the $n$-th \rep, by $X_i \super{k}$ the spatial displacement between consecutive renewals, and
 by $\tau_i \super{k}$ the corresponding temporal displacement, i.e.
 $$ S_n \super{k} \ = \ S_0 \super{k} \ + \ \sum_{i = 1}^n \ X_i \super{k}  \qquad \textrm{and} \qquad
    T_n \super{k} \ = \ \theta_k \ + \ \sum_{i = 1}^n \ \tau_i \super{k}. $$

\noindent
 The next step is to show that for $t$ large enough we can find arbitrarily many sites occupied by the particles $\gamma_t^x (k)$, $k \geq 1$, i.e. there are
 infinitely many particles that do never coalesce.
 Before giving the idea of the proof, we start by proving some preliminary results intending to describe the behaviour of the particles.
 We give the proof for the \fpa \ only but the same holds for the other members of the system.

\begin{lemma}
\label{drift}
 Let $K$ be a positive integer and $B_K$ the spatial box $[-K, K]^d$.
 Then for all $t \geq \theta_1$ and $z \in \Z^d$ the event $A$ = \{ the \fpa \ leaves $z + B_K$ in less than one unit of time and then survives without giving
 birth until $\,t + K$ \} has positive probability.
\end{lemma}

\begin{proof}
 If the particle stands out of $z + B_K$ at time $t$, then $P(A)$ can be bounded from below by $e^{-(\lambda + 1) K}$ that is the probability that the particle
 lives without giving birth for $K$ units of time.
 Otherwise, we require the \fpa \ to increase or decrease its first spatial coordinate until leaving $z + B_K$.
 Taking $\delta > 0$ such as $(K + 2) \,\delta < 1$, a straighforward calculation shows that
 $$ P \,(A) \ \geq \
    e^{-(1 + \lambda) K} \ [\,e^{-2 \delta} \ (1 - e^{- \delta}) \ e^{- \lambda \delta} \ |\mathcal V|^{-1} \ (1 - e^{- \lambda \delta}) \,]^K \ > \ 0 $$
 since it takes at most $K$ steps for the particle to reach a face.
 This concludes the proof.
\end{proof}

\begin{lemma}
\label{cubes}
 For $t \geq \theta_1$ and $K > 0$, denote by $H_{t, K}$ the event that the triangles of the \fpa \ are contained in some box $z + B_K$ from time $t$
 and for $K$ units of time.
 Then, for any $\varepsilon > 0$, there exists a large enough $K$ so that $H_{t, K}$ occurs with probability at least $\,1 - \varepsilon$.
\end{lemma}

\begin{proof}
 To make the notations easier, we will omit, all along the proof, the superscripts 1 that refer to the \fpa.
 Moreover, since the graphical representation is translation invariant, we can suppose that the first \rep \ after time $t$ stands on the origin and let $z = 0$.
 Observe that in this case, $H_{t, K}$ occurs if the spatial locations of the particle at the renewals do not leave the box $\frac{1}{2} \,B_K$
 and if each triangle is smaller than $\frac{1}{2} \,K$.
 Denoting by $\hat \gamma_s$ the location of the random walk $S_n$ at time $s$ and by $N_{t, K}$ the number of triangles between time $t$ and $t + K$,
 we then obtain
\begin{align*}
 & P \,(\,H_{t, K} \textrm{ does not occur}\,) \ \leq \ P \,\bigg(\,N_{t, K} > \frac{2}{m} \,K \,\bigg) \\
 & \hspace{05mm} + \ P \,\bigg(\,\textrm{the largest triangle is } > \frac{K}{2} \ ; \ N_{t, K} \leq \frac{2}{m} \,K \,\bigg) \\
 & \hspace{10mm} + \ P \,\bigg(\,|\!|\hat \gamma_s|\!|_{_{\infty}} > \frac{K}{2} \textrm{ for some } \ t \leq s \leq t + K \ ;
   \ N_{t, K} \leq \frac{2}{m} \,K \,\bigg)
\end{align*}
 where $m = \E \,\tau_1$.
 To begin with, we can bound the first term on the right-hand side by using the large deviation estimate
 $$ P \,\bigg(\,N_{t, K} > \frac{2}{m} \,K \,\bigg) \ \leq \ C \,e^{- \beta \,K}. $$
 Moreover, since the random vectors $X_i$, $i \geq 1$, have the same distribution, the second term can be bounded as follows
\begin{align*}
 & P \,\bigg(\,\textrm{the largest triangle is } > \frac{K}{2} \ ; \ N_{t, K} \leq \frac{2}{m} \,K \,\bigg) \\
 & \hspace{10mm} \leq \ P \,\bigg(\max_{1 \leq i \leq 2m^{-1}K} |\!|X_i|\!|_{_{\infty}} > \frac{K}{2} \,\bigg) \
   \leq \ \frac{2}{m} \,K \ P \,\bigg(\,|\!|X_1|\!|_{_{\infty}} > \frac{K}{2} \,\bigg) \
   \leq \ \frac{2}{m} \,K \,C \,e^{- \frac{1}{2} \,\beta \,K}.
\end{align*}
 Observe now that if the random walk $S_n$ leaves the box $\frac{1}{2} \,B_K$ then at least one of its coordinates is bigger than $\frac{1}{2} \,K$ so
 the third term can be bounded by
\begin{align*}
 & \sum_{i = 1}^d \ P \,\bigg(\,|\hat \gamma_s\super{(i)}| > \frac{K}{2} \textrm{ for some } t \leq s \leq t + K \ ;
   \ N_{t, K} \leq \frac{2}{m} \,K \,\bigg) \\
 & \hspace{10mm} = \ d \ P \,\bigg(\,|\hat \gamma_s\super{(1)}| > \frac{K}{2} \textrm{ for some } t \leq s \leq t + K \ ;
   \ N_{t, K} \leq \frac{2}{m} \,K \,\bigg)
\end{align*}
 where the superscript $i$ refers to the $i$-th coordinate.
 We can use the reflection principle to bound this last term by
 $$ 2 \,d \ P \,\bigg(\,|S_{2m^{-1}K}\super{(1)} - S_0\super{(1)}| > \frac{K}{2} \,\bigg). $$
 On the other hand, Chebyshev's inequality gives for any $\theta > 0$
 $$ P \,\bigg(\,|S_{2m^{-1}K}\super{(1)} - S_0\super{(1)}| > \frac{K}{2} \,\bigg) \ \leq \
    e^{- \frac{1}{2} \,\theta \,K} \ \prod_{i = 1}^{2m^{-1}K} \,\E \,e^{\,\theta \,X_i\super{(1)}} \ = \
    e^{- \frac{1}{2} \,\theta \,K} \cdot e^{\,2m^{-1}K \,\log \phi(\theta)} $$
 where $\phi(\theta)$ is the moment generating function of $X_1\super{(1)}$.
 Since $\E \,X_1\super{(1)} = 0$ and $\Var \,X_1\super{(1)} \! < \infty$ we can state that $\log \phi(\theta) \leq C \,\theta^2$ for some $C > 0$ and
 for $\theta$ small enough.
 In particular, taking $\theta = \displaystyle\frac{1}{\sqrt K}$ in the last expression we conclude that
 $$ P \,\bigg(\,|\!|\hat \gamma_s|\!|_{_{\infty}} > \frac{K}{2} \textrm{ for some } \ t \leq s \leq t + K \ ; \ N_{t, K} \leq \frac{2}{m} \,K \,\bigg) \
    \leq \ 2 \,d \,e^{- \frac{1}{4} \,\sqrt K} $$
 for $K$ large enough.
 Putting things together, we can finally maintain that
 $$ P \,(H_{t, K}) \ \geq \ 1 \ - \ C \,e^{- \beta \,K} \ - \ \frac{2}{m} \,K \,C \,e^{- \frac{1}{2} \,\beta \,K} \ - \
    2 \,d \ e^{- \frac{1}{4} \sqrt K}. $$
 This completes the proof of the lemma.
\end{proof}

\noindent
 We now prove that for $t$ large enough the number of sites occupied at time $t$ by the particles $\gamma_t^x (k)$ can be made arbitrarily large.
 The strategy of the proof consists in extracting a subsequence of ancestors that do never coalesce.
 To do this, we proceed by induction.
 Let $n \geq 1$ be a positive integer and suppose that the $n$ particles $\gamma_t^x (k_i)$, $1 \leq i \leq n$, never collide.
 Our main ingredient to prove that there exists a.s. another particle $\gamma_t^x (\ell)$ that does not coalesce with the first $n$ ones is given by the
 following lemma.

\begin{lemma}
\label{coalescence}
 For any $\varepsilon > 0$, there exists a large enough $K_n$ such as the following holds:
 If for all $1 \leq i \leq n$ both particles $\gamma_t^x (k_i)$ and $\gamma_t^x (\ell)$ are good and get separated by a distance $\geq K_n$ at some
 time $t \geq 0$ then
 $$ P \,(\,\textit{the particles do never coalesce}\,) \ \geq \ 1 - \varepsilon. $$
 Moreover, on the set where coalescence does not occur, the distance between the particles tends to infinity as $t \rightarrow \infty$.
\end{lemma}

\begin{proof}
 Let $1 \leq i \leq n$ and suppose that the particles $\gamma_t^x (k_i)$ and $\gamma_t^x (\ell)$ are good and separated by a distance $\geq K$ at
 some time $t \geq 0$.
 Then, Lemma 5.5 of Neuhauser (1992) implies that there exists a constant $C > 0$ such as
\begin{align*}
 & P \,\Big(\,\textrm{the particles } \gamma_t^x (k_i) \textrm{ and } \gamma_t^x (\ell) \textrm{ do never coalesce}\,\Big) \\
 & \hspace{10mm} = \ P \,\Big(\,\lim_{t \rightarrow \infty} \,|\!|\gamma_t^x (k_i) - \gamma_t^x (\ell)|\!| \, = \, + \,\infty \,\Big) \
   \geq \ 1 \ - \ C \,K\super{- 1 / 10} - \ 2 \,C \,K\super{- 3 / 32}.
\end{align*}
 In particular, for $K_n$ large enough,
\begin{align*}
 & P \,\Big(\,\gamma_t^x (k_i) \textrm{ and } \gamma_t^x (\ell) \textrm{ coalesce for some } 1 \leq i \leq n \,\Big) \\
 & \hspace{10mm} \leq \ \sum_{i = 1}^n \ P \,\Big(\,\gamma_t^x (k_i) \textrm{ and } \gamma_t^x (\ell) \textrm{ coalesce} \,\Big) \
   \leq \ n \,C \,K_n\super{- 1 / 10} + \ 2 \,n \,C \,K_n\super{- 3 / 32} \leq \ \varepsilon.
\end{align*}
 This proves the lemma.
\end{proof}

\noindent
 Putting the last three lemmas together, we can prove the following

\begin{lemma}
 Let $A_{\ell, K}$ = \{ for every $1 \leq i \leq n$, both particles $\gamma_t^x (k_i)$ and $\gamma_t^x (\ell)$ are good and get separated by a distance
 $\geq K$ at least once between time $\theta_{\ell}$ and $\theta_{\ell} + K$ \}.
 Then for $K$ large enough, the event $A_{\ell, K}$ occurs for infinitely many $\ell \geq 1$.
\end{lemma}

\begin{proof}
 Since the first $n$ particles do not coalesce and that $\lim \theta_{\ell} = + \infty$, we can find by Lemma \ref{coalescence} a large enough $\ell$
 such as the distances between the particles at time $\theta_{\ell}$ are bigger than $4 \sqrt d K$.
 Then, denote by $H_i$ the event that the triangles of the $k_i$-th particle are contained, between time $\theta_{\ell}$ and $\theta_{\ell} + K$, in
 some box $\Omega_i = z_i + B_K.$
 Observe that $\ell$ has been chosen so that $\Omega_i \cap \Omega_j = \varnothing$ as soon as $i \neq j$.
 Hence the events $H_i$ are determined by disjoint parts of the graph, and then are independent.
 This together with Lemma \ref{cubes} implies that for $K$ large enough
 $$ P \,(\,H_1 \,\cap \,H_2 \,\cap \,\cdots \,\cap \,H_n \,) \ \geq \ (1 - \varepsilon)^n. $$
 Now that the first $n$ particles are trapped inside large disjoint cubes, we require each of them to be good at least once between time $\theta_{\ell}$
 and $\theta_{\ell} + K$.
 Since this occurs if the $n$ particles have at least one renewal in this time laps, the probability of this event can be bounded from below by
 $$ \prod_{i = 1}^n \,P \,(\,\tau_1\super{k_i} < K - 1 \,) \ \geq \ \Big[\,1 - C \,\exp \,(- \beta \,(K - 1))\,\Big]^n $$
 where $C$ and $\beta$ are some positive constants.
 The last thing we need is that the $\ell$-th particle is good and its distance from each other ancestor is $\geq K$ between $\theta_{\ell} + 1$
 and $\theta_{\ell} + K$.
 Since the cubes $\Omega_i$ are separated by a distance $\geq 2K$, this obviously occurs, by Lemma \ref{drift}, with positive probability.
 Putting things together we can state that there exists an $\eta > 0$ such as $P(A_{\ell, K}) \geq \eta$ for $\ell$ large enough.
 Observing finally that the events $A_{\ell_1, K}$ and $A_{\ell_2, K}$ are independent as soon \linebreak as $|\theta_{\ell_1} - \theta_{\ell_2}| > K$,
 we can conclude by the Borel-Cantelli Lemma that a.s. the events $A_{\ell, K}$ occur for infinitely many $\ell \geq 1$.
\end{proof}

\noindent
 On the event $A_{\ell, K}$ coalescence occurs with probability at most $\varepsilon$ for $K$ large enough.
 Moreover, $A_{\ell, K}$ occurs for infinitely many $\ell \geq 1$, so we eventually find a particle $\gamma_t^x (\ell)$ that never coalesces with the first $n$
 members of the sequence.
 Since $n$ has been chosen arbitrarily, this proves that there are infinitely many particles that do never coalesce.
 To conclude the proof of Theorem \ref{diagonale}, we use again the dual process $\tilde \eta_s\super{(x, t)}$, $0 \leq s \leq t$, starting at $(x, t)$.
 Since, by choosing $t$ large enough, we can find arbitrarily many different sites occupied by the particles $\gamma_t^x (k)$ and that $\eta_0$ is translation
 invariant, there exists a.s. a lower ancestor $\gamma_t^x (n)$ that lands, at time 0, on a red site.
 Now, look at the ancestors that come before $\gamma_t^x (2)$ in the hierarchy until we find a particle that lands on an occupied site.
 If it is a red site, it will bring a \gpa \ to $(x, t)$ and the proof is done.
 If it is a green site, $S_0\super{1}$ will be vacant at time $(t - \theta_1)\super{+}$ whatever the color of the next ancestors.
 Then, we look at $\gamma_t^x (2)$, and so on.
 If however not one of the ancestors that come before $\gamma_t^x (n)$ succeeds in painting $(x, t)$ in green, this last one will do it.
 This completes the proof of Theorem \ref{diagonale}.


\section{\normalsize\bf{Proof of Theorem \ref{haut}}}
\label{coexistence}

\noindent
 This last section is devoted to the proof of Theorem \ref{haut}.
 In particular, we will prove that, in dimension 2, coexistence occurs
 for an open set of values $(\lambda_1, \lambda_2)$ in $\R^2$.
 We conjecture that such a property holds in any dimension but our proof heavily relies on results of Durrett and Neuhauser (1991) that have been proved
 in $d = 2$ only.
 Before going into the details of the proof, we start by describing the setting of their results.
 First of all, let $\lambda_2 = 0$.
 We recall that in this context, 0 can be interpreted as a living tree, 1 as a burning tree and 2 as a burnt site.
 See Sect. \ref{introduction} for more details.
 We start the process with one burning tree in a virgin forest and wonder whether the fire will go out or not.
 To do this, let $B = (-L, L)^2$ and, for any $m \in \Z$, $B_m = m L e_1 + B$, where $e_1 = (1, 0)$ is the first unit vector.
 Then, denote by $\mathcal L$ the subset of $\Z^2$ given by
 $$ \mathcal L \ = \ \{\,(m, n) \in \Z^2 \,,\, m + n \textrm{ is even} \,\} $$
 and say that a site $(m, n) \in \mathcal L$ is \emph{occupied} if the following two conditions are satisfied.
\begin{enumerate}
\item [(1)] There are more than $L\super{1 / 2}$ burning trees in $B_m$ at some time $t \in [n \,\Gamma \,L, (n + 1) \,\Gamma \,L]$.
\item [(2)] There is at least one burning tree in $B_m$ at all times $t \in [(n + 1) \,\Gamma \,L, (n + 2) \,\Gamma \,L]$.
\end{enumerate}
 Here $\Gamma$ is a positive constant that will be fixed later.
 The following lemma is a comparison result of Durrett and Neuhauser (1991).
 It guarantees, in the setting of the forest fire model, that coexistence occurs for $\lambda_1$ greater than some critical value
 $\alpha_c \in (0, \infty)$.

\begin{lemma}
\label{fire}
 If $\lambda_2 = 0$ and $\lambda_1 > \alpha_c$ then $\Gamma$ and $L$ can be chosen so that the set of occupied sites dominates the set of wet sites in
 a one-dependent oriented percolation process on $\mathcal L$ with parameter $p = 1 - 6^{\,- 36}$.
\end{lemma}

\noindent
 The next step is to show that the conclusion of Lemma \ref{fire} holds for small values of the parameter $\lambda_2$.
 To do this, take a site $(m, n) \in \mathcal L$, and denote by $B_{m, n}\super{+}$ the space-time region
 $$ B_{m, n}\super{+} \ = \ B_m \,\times \,[n \,\Gamma \,L, (n + 2) \,\Gamma \,L] \ \cup \ B_{m + 1} \,\times \,[(n + 1) \,\Gamma \,L, (n + 3) \,\Gamma \,L]. $$
 Let $H_{m, n}\super{+}$ be the event that the \gpa s do not give birth in the box $B_{m, n}\super{+}$.
 A straightforward calculation implies that the volume of $B_{m, n}\super{+}$ is bounded by $14 \,\Gamma \,L^3$ so that
 $$ P \,(H_{m, n}\super{+}) \ \geq \ e^{- 14 \,\lambda_2 \,\Gamma \,L^3} \ \geq \ 1 - 6^{\,- 36} $$
 for $\lambda_2 > 0$ small enough.
 Denoting by $B_{m, n}\super{-}$ the set
 $$ B_{m, n}\super{-} \ = \ B_m \,\times \,[n \,\Gamma \,L, (n + 2) \,\Gamma \,L] \ \cup \ B_{m - 1} \,\times \,[(n + 1) \,\Gamma \,L, (n + 3) \,\Gamma \,L] $$
 and by $H_{m, n}\super{-}$ the event that the \gpa s do not give birth in $B_{m, n}\super{-}$, the same arguments obviously imply that
 $$ P \,(H_{m, n}\super{-}) \ \geq \ 1 - 6^{\,- 36}. $$
 On the other hand, Lemma \ref{fire} states that if $(m, n)$ is occupied and $H_{m, n}\super{+}$ (resp. $H_{m, n}\super{-}$) occurs, the probability
 that $(m + 1, n + 1)$ (resp. $(m - 1, n + 1)$) is occupied can be bounded from below by $1 - 6^{\,- 36}$.
 Putting things together, we can finally maintain that, for $\lambda_2 > 0$ small enough, the conclusion of \ref{fire} holds for
 $p = 1 - 2 \times 6^{\,- 36}$. \\

\noindent
 To conclude the proof of Theorem \ref{haut}, observe first of all that classical results about percolation imply that for such a value of $p$ there
 exists a.s. an infinite cluster of wet sites in $\mathcal L$.
 See for instance Durrett (1984), Sect. 10.
 Next, start the process from an initial configuration that has at least $L\super{1 / 2}$ \rpa s in each box $B_m$ and extract a convergent subsequence of
 the Cesaro averages.
 The limit $\mu$ we then obtain is known to be a stationary distribution.
 Moreover, by comparison with the percolation process defined above, $\mu$ concentrates on configurations with infinitely many \rpa s.
 This concludes the proof of Theorem \ref{haut}. \\



\noindent
\textbf{Acknoledgement}.
 I would like to thank my supervisors Claudio Landim who, as usual, suggests very interesting research subjects,
 and Olivier Benois for his encouragement, his advice, and the time he has devoted to me.
 I am also grateful to Pierre Margerie, from the Laboratoire d'\'Ecologie de l'Universit\'e de Rouen, groupe ECODIV, for his precious explanations
 relating to ecological phenomena. \\


\vspace{20pt}

\hfill\begin{minipage}{170pt}
\small{\textsc{D\'epartement de Math\'ematiques \\
 Universit\'e de Rouen \\
 76128 Mont Saint Aignan, France \\
 E-mail}: nicolas.lanchier@univ-rouen.fr}
\end{minipage}

\end{document}